\documentclass[a4paper,11pt]{amsart}
\usepackage{epsfig}
\usepackage{amsmath, amsfonts }
\usepackage{amssymb}
\usepackage{amscd}
\usepackage{eufrak}
\usepackage{amssymb}
\address{\'Ecole Nationale Sup\'erieure Polytechnique,\\
B.P.~8390 Yaound\'e, Cameroun\footnote{currently at UMR 5030 (CNRS), D\'epartement des Sciences Math\'ematiques Universit\'e Montpellier II Case courrier 051-Place Eug\'ene Bataillon 34095 Montpellier CEDEX 05, France e-mail:tbouetou@darboux.math.univ-montp2.fr}}
\email{tbouetou@polytech.uninet.cm}

\title{CLASSIFICATION OF SOLVABLE 3-DIMENSIONAL LIE TRIPLE SYSTEMS}

\author{Thomas B. Bouetou}
\begin{document}
\maketitle
\begin{abstract}
We give the classification of solvable and splitting Lie triple systems and it turn that, up to isomorphism there exist 7 non isomorphic canonical Lie triple systems and, 6 non isomorphic splitting canonical Lie triple systems and find the solvable Lie algebras associated.
\end{abstract}


\section{INTRODUCTION}

\em A Lie triple system(LTS), is a space where is defined a ternary operation, verifying some conditions, namely the Jacobi identity and the derivation identity. 
They where first introduce by Jacobson \cite{ja}. Later on Lister \cite{lis}
 gave a structure theory and the classification of simple LTS. 
Yamaguty \cite{ya} obtain from a total geodesic space triple algebras which
where the generalization of LTS. Loos \cite{lo} show that a symmetric space can
be seen as a quasigroup, and Sabinin \cite{sab1,sab2} show that a quasigroup can be
seen as a homogeneous space. In particular, any Bol loop under the left action
 derivative give a LTS i.e. the description of the infinitesimal structure
 of a smooth Bol loop contain a LTS. This fact give the idea of investigation
 of LTS since also the use of  LTS appear in the ordinary differential equation
 functional analysis...In this paper our main object is to give the 
classification of solvable and splitting LTS up to isomorphism our approach is 
based on  the enveloping Lie algebras of a LTS.Since the Lie algebras obtain 
from the standard embedding of a LTS it is an enveloping Lie algebra i.e. if a
LTS is solvable his enveloping Lie algebra is solvable, conversely if a Lie
 algebra is solvable the the LTS obtain is solvable. Considering the 
classification of solvable Lie algebra, we will carry out the classification 
of LTS of small dimension.

We will organize this pepper as follows: The first part is the introduction,
the second part we give the definition and some result about LTS. In the third
 part we give the classification of LTS of dimension two, the forth part we 
give the classification of solvable LTS and finally the last part we give the 
classification of splitting LTS.

\section{ABOUT LIE TRIPLE SYSTEMS}

\underline{{\bf Definition 2.1}} The vector space $ \mathfrak{M} $ 
(finite over the field of real numbers $ \mathbb{R}$) with trilinear operation (x,y,z) is called a LTS  if the following identities are verify:
$$
(x,x,y)=0
$$
$$
(x,y,z)+(y,z,x)+(z,x,y)=0
$$
$$
(x,y,(u,v,w))=((x,y,u),v,w)+(u,(x,y,v),w)+(u,v,(x,y,w))
$$
Let $ \mathfrak{M} $ be a  LTS, a subspace
 $ \mathfrak{D}  \subset  \mathfrak{M} $ is called a subsystem if
 $ (\mathfrak{D}, \mathfrak{D}, \mathfrak{D})\subset \mathfrak{D} $, and is
called an
 ideal, if $ ( \mathfrak{D}, \mathfrak{M}, \mathfrak{M}) \subset \mathfrak{D} $.
 The ideals are the Kernel of the homomorphism of the  LTS see \cite{lo, trof}.

\underline{\bf{Example}} For a typical way of construction of a  LTS see
 in  \cite{lo,trof}.

Let $ \mathfrak{G}$ be a Lie algebra (finite over the field of real numbers
$ \mathbb{R}) $  and $ \sigma $-an involutive automorphism, then
$$
\mathfrak{G}=\mathfrak{G}^{+} \dotplus \mathfrak{G}^{-}
$$
where $ \sigma | \mathfrak{G^{+}}=Id $ and  $ \sigma | \mathfrak{G^{-}}=-Id $,
as  any element $x$ from $ \mathfrak{G}$ can be written in the form:
$$
x=\frac{1}{2}(x+\sigma x)+\frac{1}{2}(x-\sigma x),
$$
where  $x+\sigma x  \in \mathfrak{G}^{+} $,  $ x-\sigma x \in \mathfrak{G}^{-}$ and $ \mathfrak{G}^{+} \cap \mathfrak{G}^{-}=0 $.

The following inclusions hold:

$$
[\mathfrak{G^{+}}, \mathfrak{G^{+}}] \subset \mathfrak{G^{+}}, [\mathfrak{G^{+}}, \mathfrak{G^{-}}] \subset \mathfrak{G^{-}}, [\mathfrak{G^{-}}, \mathfrak{G^{-}}] \subset \mathfrak{G^{+}}.
$$
 Then the subspace $ \mathfrak{G^{-}} $ turns into a  LTS 
relatively under the operation $ (x,y,z)=[[x,y],z] $.

The inverse construction \cite{lo}.

Let $\mathfrak{M}$ be a  LTS and define by 
$$h(X,Y): z \longrightarrow (X,Y,Z)$$
a linear transformation of the space $\mathfrak{M}$ into itself where 
$X,Y,Z \in \mathfrak{M}$.

Let $H$ be a subspace of the space of linear transformations of the  LTS  $\mathfrak{M}$ whose elements are the transformations of the form $h(X,Y)$. The 
vector space $\mathfrak{G}=\mathfrak{M} \dotplus H$, become a Lie algebra 
relatively to the commutator $[A,B]=AB-BA$, $[A,X]=-[X,A]=AX$;
$[X,Y]=h(X,Y)$ where $A,B \in H$, $X,Y \in \mathfrak{M}$.

Let us define the mapping $\sigma$ with the condition
$\sigma(A)=A,$ if $ A\in H$ and $\sigma(X)=-X, X \in \mathfrak{M}$, then
$\sigma$ is an involutive automorphism of a Lie algebra $\mathfrak{G}=\mathfrak{M} \dotplus H$.

The algebra $\mathfrak{G}$ constructed above from the  LTS, is
called universal enveloping Lie algebra of the  LTS $\mathfrak{M}$.

\underline{\bf Definition 2.2} The derivation of the  LTS 
$\mathfrak{M}$, is called
 the linear transformation
$\mathfrak{d}: \mathfrak{M} \longrightarrow \mathfrak{M}$ such that
$$ (X,Y,Z)\mathfrak{d}=(X\mathfrak{d},Y,Z)+(X,Y\mathfrak{d},Z)+(X,Y,Z\mathfrak{d}).$$

One can verify that, the set $\mathfrak{d}(\mathfrak{M})$ of all the 
derivation of the  LTS  $\mathfrak{M}$ is a Lie algebra of the 
linear transformations acting on $\mathfrak{M}.$

\underline{\bf Definition 2.3} The embedding of a  LTS 
$\mathfrak{M}$ into a Lie 
algebra $\mathfrak{G}$ is called the linear injection
$R: \mathfrak{M} \longrightarrow \mathfrak{G}$ such that 
$(X,Y,Z)=[[X^R,Y^R],Z^R]$.

The embedding $R$ of the  LTS  $\mathfrak{M}$ into the Lie algebra
$\mathfrak{G}$ is called canonical, if the envelope of the image of the set 
$\mathfrak{M}^R$ in the Lie algebra $\mathfrak{G}$ coincide with 
$\mathfrak{G}$ and $h$ does not contain trivial ideals of Lie algebra
$\mathfrak{G}.$ Let us note that if the  LTS  $\mathfrak{M}$ is
a subset of the Lie algebra $\mathfrak{G}$, then $(X,Y,Z)=[[X,Y],Z]$
and $[\mathfrak{M},\mathfrak{M}]$ is a subalgebra of the Lie algebra 
$\mathfrak{G}$ hence
$\mathfrak{M}+[\mathfrak{M},\mathfrak{M}]$- is a Lie subalgebra of 
$\mathfrak{G}$ and the initial embedding $R$ can be consider as canonical in 
$\mathfrak{M}^R+[\mathfrak{M}^R,\mathfrak{M}^R]$; this lead us to formulate
the following proposition:

\underline{{\bf Proposition 2.1}} For any finite  LTS 
$\mathfrak{M}$
over $\mathbb{R}$, there exist one and only up to automorphism accuracy, one 
canonical embedding to the Lie algebra.

\subsection{SOLVABLE AND SEMISIMPLE LIE TRIPLE SYSTEM}

Following \cite{lis}:, let $\Omega$- be an ideal of the  LTS 
$\mathfrak{M}$, we assume $\Omega^{(1)}=(\mathfrak{M},\Omega,\Omega)$ and,
  $\Omega^{(k)}=(\mathfrak{M},\Omega^{(k-1)},\Omega^{(k-1)})$

\underline{{\bf Proposition 2.2}} \cite{lis} For all natural number $k$, the subspace
$\Omega^{(k)}$ is an ideal of $\mathfrak{M}$ and we have the following
inclusions:
$$ \Omega \supseteq \Omega^{(1)} \supseteq ........\supseteq \Omega^{(k)}$$

Proof

$$
(\Omega^{(1)},\mathfrak{M},\mathfrak{M})=((\mathfrak{M},\Omega,\Omega),\mathfrak{M},\mathfrak{M})\subseteq ((\mathfrak{M},\Omega,\mathfrak{M}),\Omega,\mathfrak{M})+[[[\mathfrak{M},\Omega],[\mathfrak,\Omega]],\mathfrak{M}]
$$
according to the definition of a  LTS
$$
(\Omega^{(1)},\mathfrak{M},\mathfrak{M})\subseteq (\Omega,\Omega,\mathfrak{M})+[(\mathfrak{M},\Omega,\mathfrak{M}),[\mathfrak{M},\Omega]]\subseteq (\mathfrak{M},\Omega,\Omega)+(\mathfrak{M},\Omega,\Omega)=\Omega^{(1)}
$$
that means  $ \Omega^{(1)}$ is an ideal of $\mathfrak{M}$
further more $\Omega^{(k)}=(\Omega^{(k-1)})^{(1)}$ hence each $\Omega^{(i)}$ is 
an 
ideal in $\mathfrak{M}$.

\underline{{\bf Definition 2.4}} The ideal $\Omega$ of a  LTS  
$\mathfrak{M}$ is
called solvable, if there exist a natural number $k$  such that $\Omega^{(k)}=0.$

\underline{{\bf Proposition 2.3}} \cite{lis} If $\Omega$ and $\Theta$ are two solvable ideals
 of a  LTS  $\mathfrak{M}$ then  $\Omega+\Theta$ is also a solvable
ideal in $\mathfrak{M}$.

Proof

using the definition of a  LTS, the following inclusion hold:
$(\Theta+\Omega)^{(1)}\subseteq (\mathfrak{M},\Theta,\Theta)+(\mathfrak{M},\Omega,\Omega)+(\mathfrak{M},\Theta,\Omega)+(\mathfrak{M},\Omega,\Theta)\subseteq \Theta^{(1)}+\Omega^{(1)}+\Theta \cap \Omega$.
 
Assume for every natural number $k$ the following inclusion holds:

$(\Theta+\Omega)^{(k)}\subseteq \Theta^{(k)}+\Omega^{(k)}+\Theta \cap \Omega$

by induction let's prove that its holds for $(k+1)$

$(\Theta+\Omega)^{(k+1)}=(\mathfrak{M},(\Theta+\Omega)^{(k)},(\Theta+\Omega)^{(k)})\subseteq(\mathfrak{M},\Theta^{(k)}+\Omega^{(k)}+\Theta \cap \Omega,(\Theta \cap \Omega))\subseteq \Theta^{(k+1)}+\Omega^{(k+1)}+\Theta \cap \Omega$

hence the result

\underline{{\bf Definition 2.5}} The radical of a  LTS  denoted 
by 
$\mathfrak{R}(\mathfrak{M})$, is called the maximal solvable ideal of the  LTS  $\mathfrak{M}.$

A  LTS  $\mathfrak{M}$ is called semi-simple if $\mathfrak{R}(\mathfrak{M})=0$.

\underline{{\bf Theorem 2.1}} \cite{lis} If $\mathfrak{R}$ is a radical in $\mathfrak{M}$ 
then $(\mathfrak{M}\setminus \mathfrak{R})$ is semisimple. And if 
$\Omega $ is an ideal in $\mathfrak{M}$ such that $(\mathfrak{M}\setminus \mathfrak{R})$ is semisimple then $\Omega \supset \mathfrak{R}$.   

\underline{{\bf Proposition 2.4}} \cite{lis} The enveloping Lie algebra, of 
a solvable  LTS  is solvable. And if a  LTS  has some solvable
enveloping Lie algebra, it is solvable.

\underline{{\bf Theorem 2.2}} If $\mathfrak{M}$ is a semisimple  LTS, then
the universal enveloping Lie algebra $\mathfrak{G}$ is semisimple.

\underline{{\bf Theorem 2.3}} \cite{bou} Let $\mathfrak{M}$ be a  LTS  and
$\mathfrak{G}=\mathfrak{M} \dotplus \mathfrak{h}$ his canonical enveloping Lie 
algebra and $\mathfrak{r}$- the radical of the Lie algebra $\mathfrak{G}$.
In $\mathfrak{G}$ there exist a subalgebra $\mathfrak{P}$ semisimple 
supplementary to with $\mathfrak{r}$ such that:
$$
\mathfrak{M}=\mathfrak{M}' \dotplus \mathfrak{M}''\; (direct\; sum\; of\; vectors \;spaces)
$$
where
$$
\mathfrak{M}'=\mathfrak{M} \cap \mathfrak{r} -radical\; of\; the\;  LTS\; \mathfrak{M}
$$
$$
\mathfrak{M}''=\mathfrak{M} \cap \mathfrak{P} -semisimple\; subalgebra\; of\;  LTS \; \mathfrak{M}
$$
$$
\mathfrak{h}=\mathfrak{h}' \dotplus \mathfrak{h}''\; (direct\; sum\; of\; vectors \; spaces)
$$
$$
\mathfrak{h}'=\mathfrak{h} \cap \mathfrak{r}
$$
and
$$
\mathfrak{h}''=\mathfrak{h} \cap \mathfrak{P}\; are\; subalgebra \;in \; \mathfrak{h}
$$
$$
\mathfrak{r}=\mathfrak{M}' \dotplus \mathfrak{h}'
$$
$$
\mathfrak{P}=\mathfrak{M}'' \dotplus \mathfrak{h}''.
$$

\subsection{PROBLEM SETTING}

Let $\mathfrak{M}$ be a  LTS  and $dim \mathfrak{M}=3$. To be
consistent with the above Theorem the following cases are possible:
\begin{enumerate}
\item  semisimple case\\
$\mathfrak{M}$- semisimple  LTS (in fact simple).
 About the classification of such  LTS see \cite{ber,fed1,lis}

\item Solvable case\\
$\mathfrak{M}$ is a solvable  LTS. The classification of such 
system is given  section 4. 
\item Splitting case\\
$$\mathfrak{M}=\mathfrak{M}_1 \dotplus \mathfrak{M}_2 $$
where $\mathfrak{M}\equiv \mathbb{R}$- solvable ideal of dimension 1 in $\mathbb{R}$
and $\mathfrak{M}_2$ -semisimple  LTS of dimension 2 This type of  LTS is  considered at the last section.
\end{enumerate}

\section{CLASSIFICATION OF LIE TRIPLE SYSTEM OF DIMENSION 2}

For a better survey of such  LTS, we will write their trilimear
operation in a special form.

Let $\mathfrak{M}$ be a 2-dimensional  LTS we write the trilinear
operation $(X,Y,Z)=\beta (X,Y)Y-\beta (Y,Z)X$
where $\beta:V \times V\longrightarrow \mathbb{R}$ is a symmetric form. The
choice of the basis $V=<e_1,e_2>$ one can reduce the symmetric form to the
view: 
\begin{displaymath}
\left(\begin{array}{cc}
\alpha & 0 \\
0 & \nu \\
\end{array}\right), 
\end{displaymath}

where $ \alpha, \nu= \pm 1;0$.

By introducing the notation of the derivation
$$
\mathfrak{D}_{x,y}:\mathfrak{M} \longrightarrow \mathfrak{M}
$$
$$
z\longmapsto (x,y,z)
$$
$$
\mathfrak{h}=\left\{\mathfrak{D}_{x,y}\right\}_{x,y \in \mathfrak{M}}.
$$
And

$\mathfrak{G}=\mathfrak{M} \dotplus \mathfrak{h}$- canonical enveloping 
Lie algebra of the  LTS $\mathfrak{M}$.

Let $\mathfrak{M}=<e_1,e_2>$ then,

$ \mathfrak{h}=\left\{tD_{x,y}\right\}_{t \in \mathbb{R}}$,
$$
e_{1}D=(e_1,e_2,e_1)=\beta (e_1,e_1)e_2
$$
$$
e_{2}D=(e_1,e_2,e_2)=-\beta (e_2,e_2)e_1
$$

$\mathfrak{G}=<e_1,e_2,e_3>$ 

where 
$[e_1,e_2]=e_3, \; [e_1,e_3]=-e_{1}D, \; [e_2,e_3]=-e_{2}D$

Therefore we can have the up to isomorphism accuracy the following five
cases:

\begin{enumerate}
\item (Spherical Geometry)\\

\begin{displaymath}
\left(\begin{array}{cc}
1 & 0 \\
0 & 1 \\
\end{array}\right), 
\end{displaymath}

$\mathfrak{G}/\mathfrak{h}\cong so(3)/so(2)$
\item (Lobatchevski Geometry)\\

\begin{displaymath}
\left(\begin{array}{cc}
-1 & 0 \\
0 & -1 \\
\end{array}\right), 
\end{displaymath}

$\mathfrak{G}/\mathfrak{h}\cong sl(2,\mathbb{R})/so(2)$
\item  LTS with non compact subalgebra $\mathfrak{h}$\\
\begin{displaymath}
\left(\begin{array}{cc}
1 & 0 \\
0 & -1 \\
\end{array}\right), 
\end{displaymath}

$\mathfrak{G}/\mathfrak{h}\cong sl(2,\mathbb{R})/\mathbb{R}$
\item Solvable case\\
\begin{itemize}
\item a) 

\begin{displaymath}
\beta=\left(\begin{array}{cc}
1 & 0 \\
0 & 0 \\
\end{array}\right), 
\end{displaymath}

$e_1 \cdot e_2=e_3, \; e_1 \cdot e_3=e_2$

(This is a Lie algebra $\mathfrak{G}$ of type $g_{3,5}(p=0)$ in \cite{mub2})
\item b)
\begin{displaymath}
\beta=\left(\begin{array}{cc}
-1 & 0 \\
0 & 0 \\
\end{array}\right), 
\end{displaymath}

$e_1 \cdot e_2=e_3, \; e_1 \cdot e_3=-e_2$

(This is a Lie algebra $\mathfrak{G}$ of type $g_{3,4}(h=-1)$ in \cite{mub2})
\end{itemize}
\item Abelian case\\
$\beta=0$ $\mathfrak{G}/\mathfrak{h}\cong (\mathbb{R})^2/\left\{0\right\}$
\end{enumerate}

\section{CLASSIFICATION OF SOLVABLE LIE TRIPLE SYSTEMS OF DIMENSION 3}

Let $\mathfrak{M}$- be a solvable  LTS  of dimension 3, and 
$\mathfrak{G} \dotplus \mathfrak{h}$ its canonical enveloping Lie algebra
then $\mathfrak{G}$ is solvable in particular $\mathfrak{G}$ posses
a characteristic ideal 
$\mathfrak{G}'=[\mathfrak{G},\mathfrak{G}]\vartriangleright \mathfrak{G}$,

$\sigma \mathfrak{G}'=\mathfrak{G}'$,
$\mathfrak{G}' \cap \mathfrak{M}=\mathfrak{M}'=(\mathfrak{M},\mathfrak{M},\mathfrak{M})$ further more $\mathfrak{h} \subset \mathfrak{G}$ since
$\mathfrak{h}=[\mathfrak{M},\mathfrak{M}]$ then 
 
$\mathfrak{G}'=[\mathfrak{G},\mathfrak{G}]=\mathfrak{M}'+ \mathfrak{h}$ where
$\mathfrak{M}'\subsetneq \mathfrak{M}$ 
 
Possible situations:

\begin{enumerate}
\item $ dim \mathfrak{M'}=0 $.
 Then $ [\mathfrak{h}, \mathfrak{M}]=\mathfrak{M'}=\{O \} $, that
 means $ \mathfrak{h} \vartriangleright \mathfrak{G} $- ideal, that is why
 $ \mathfrak{h}=\{O \} $ (since $ \mathfrak{G} $- is an enveloping
 Lie algebra) and $ \mathfrak{M}=\mathbb{R} \oplus \mathbb{R} \oplus \mathbb{R} $. 
In this case, the  LTS  is Abelian and we denote it (type I).
\item $   dim \mathfrak{M'}=1 $. Choosing the base $ e_1 $, $ e_2 $, $ e_3 $ in $ \mathfrak{M} $ such that, $ \mathfrak{M'}=<e_{1}> $ and $ \mathfrak{M}=\mathfrak{M'}+<e_{2}, e_{3}> $.

We will introduce in consideration the linear transformation 
$ A, B, C: \mathfrak{M}\longrightarrow \mathfrak{M} $, define as:
\begin{displaymath}
A=(e_{1},e_{2},-)=\left(\begin{array}{ccc}
a & b & c\\
0 & 0 & 0\\
0 & 0 & 0\\
\end{array}\right), 
B=(e_{2},e_{3},-)=\left(\begin{array}{ccc}
\alpha & \beta & \gamma\\
0 & 0 & 0\\
0 & 0 & 0\\
\end{array}\right),
\end{displaymath}
\begin{displaymath}
C=(e_{3},e_{1},-)=\left(\begin{array}{ccc}
x & -\alpha -c & y\\
0 & 0 & 0\\
0 & 0 & 0\\
\end{array}\right).
\end{displaymath} 
 
And if a skew symmetric form defined as 
$ \Phi(-,-): \mathfrak{M} \times \mathfrak{M} \longrightarrow \mathbb{R} $, 
such that $ (x,y,e_{1})= \Phi(x,y)e_{1}$. 
The dimension of $ \mathfrak{M} $ is 3, that is why there exists
 $ z \in \mathfrak{M} $, $ z \neq 0 $, such that $ \Phi(-,z)=0 $. The following
cases are possible: 
\begin{itemize}
\item b.I. The skew-symmetric form $ \Phi $ is non zero and $ z $ is parallel
to $ e_1 $ $ (z \parallel e_{1}) $, then in the base $e_1 $, $ e_2 $, $ e_3 $ 
the skew-symmetric form $ \Phi $ has the corresponding matrix:

\begin{displaymath}
\left(\begin{array}{ccc}
0 & 0 & 0\\
0 & 0 & \delta\\
0 & -\delta  & 0\\
\end{array}\right), 
\end{displaymath}

where $ \delta \neq 0 $. Adjusting $ e_3 $ to $ 1 \setminus \delta e_{3} $, then
$ \Phi( e_{2},e_{3})=1 $, $ \Phi( e_{3},e_{2})=-1 $, so that $ \alpha=1 $, 
$ a=x=0$ and

\begin{displaymath}
A=(e_{1},e_{2},-)=\left(\begin{array}{ccc}
0 & b & c\\
0 & 0 & 0\\
0 & 0 & 0\\
\end{array}\right), 
B=(e_{2},e_{3},-)=\left(\begin{array}{ccc}
1 & \beta & \gamma\\
0 & 0 & 0\\
0 & 0 & 0\\
\end{array}\right),
\end{displaymath}
\begin{displaymath}
C=(e_{3},e_{1},-)=\left(\begin{array}{ccc}
0 & -1 -c & y\\
0 & 0 & 0\\
0 & 0 & 0\\
\end{array}\right).
\end{displaymath}

The verification of the defined relation of  LTS  shows that, 
with accuracy to the choice of the vector basis $ e_2 $ and $ e_3 $, it is 
possible to afford the following realization of the operators $ A $, $ B $, $ C$
 as:

\begin{displaymath}
A=0, 
B=(e_{2},e_{3},-)=\left(\begin{array}{ccc}
1 & 0 & 0\\
0 & 0 & 0\\
0 & 0 & 0\\
\end{array}\right),
\end{displaymath}
\begin{displaymath}
C=(e_{3},e_{1},-)=\left(\begin{array}{ccc}
0 &  -1 & 0\\
0 & 0 & 0\\
0 & 0 & 0\\
\end{array}\right).
\end{displaymath}  \; \; \; \; \; \; \; \; (type VII)
\item b.II. The skew-symmetric form $ \Phi $ is non zero and $ z $ is not 
parallel to $ e_1 $, let $ z=e_{2 } $, then

\begin{displaymath}
A=(e_{1},e_{2},-)=\left(\begin{array}{ccc}
0 & b & c\\
0 & 0 & 0\\
0 & 0 & 0\\
\end{array}\right), 
B=(e_{2},e_{3},-)=\left(\begin{array}{ccc}
0 & \beta & \gamma\\
0 & 0 & 0\\
0 & 0 & 0\\
\end{array}\right),
\end{displaymath}
\begin{displaymath}
C=(e_{3},e_{1},-)=\left(\begin{array}{ccc}
-1 &  -c & y\\
0 & 0 & 0\\
0 & 0 & 0\\
\end{array}\right).
\end{displaymath}

The verification of the defined relations of  LTS, show that the
 indicated case has no realization.
\item b.III. The skew-symmetric form $ \Phi $ is trivial. By completing the
 vector $ e_1 $ with the arbitrary choose vector $ e_2 $ and $ e_3 $ up to
 the base, it is possible to realize the operator $ A $, $ B $, and $ C $:

\begin{displaymath}
A=(e_{1},e_{2},-)=\left(\begin{array}{ccc}
0 & b & c\\
0 & 0 & 0\\
0 & 0 & 0\\
\end{array}\right), 
B=(e_{2},e_{3},-)=\left(\begin{array}{ccc}
0 & \beta & \gamma\\
0 & 0 & 0\\
0 & 0 & 0\\
\end{array}\right),
\end{displaymath}
\begin{displaymath}
C=(e_{3},e_{1},-)=\left(\begin{array}{ccc}
0 &  -c & y\\
0 & 0 & 0\\
0 & 0 & 0\\
\end{array}\right).
\end{displaymath}

The verification of the defined relations of  LTS, show that by
a suitable choice of basis vectors $e_2, e_3$ the following realization
of operators $A,B,C$ is possible:
\begin{itemize}
\item Abelian Type (Type above)
\item $A=C=0$,
\begin{displaymath}
B=(e_{2},e_{3},-)=\left(\begin{array}{ccc}
0 &  0 & 1\\
0 & 0 & 0\\
0 & 0 & 0\\
\end{array}\right).
\end{displaymath}\; \; \; \; \; \; \; \; \;  (Type II)

This  LTS, is obtained by a direct multiplication of a  LTS of dimension 
two $<e_1,e_2>$, by an Abelian one dimensional $<e_3>$.
\item - \begin{displaymath}
A=(e_{1},e_{2},-)=\left(\begin{array}{ccc}
0 & \pm 1 & 0\\
0 & 0 & 0\\
0 & 0 & 0\\
\end{array}\right)
\end{displaymath} B=C=0. \; \; \; \; \; \; \;  (Type III)
\item \begin{displaymath}
A=(e_{1},e_{2},-)=\left(\begin{array}{ccc}
0 &  \pm 1 & 1\\
0 & 0 & 0\\
0 & 0 & 0\\
\end{array}\right)
\end{displaymath}, B=0,
\begin{displaymath}
C=(e_{3},e_{1},-)=\left(\begin{array}{ccc}
0 &  -1 & \mp 1\\
0 & 0 & 0\\
0 & 0 & 0\\
\end{array}\right).
\end{displaymath} \; \; \; \; \; \; \; \; \; \; (Type IV)
\end{itemize}

\end{itemize} 
\item $dim \mathfrak{M}'=2$ in particular, $ \mathfrak{M}'$ is a
subsystem of dimension two in $\mathfrak{M}$. one can consider
(refer to Section 3 ) $\forall a,b,c \in \mathfrak{M}'$
$$
(a,b,c)=\beta (a,c)b-\beta (b,c)a
$$
where 
\begin{displaymath}
\beta=\left(\begin{array}{cc}
\pm 1 & 0\\
0 & 0 \\
\end{array}\right)
\end{displaymath} 
and $\mathfrak{M}'$ is a two-dimensional Abelian ideal in $\mathfrak{M}$. 
 In the first case the choice of the base $ \mathfrak{M}=<e_{1}, e_{2},e_{3}>$ 
such that $ \mathfrak{M'}=<e_{1},e_{2}> $, the operations of the  LTS are 
reduced to:

\begin{displaymath}
A=(e_{1},e_{2},-)=\left(\begin{array}{ccc}
0 & \pm 1 & x\\
0 & 0 & y\\
0 & 0 & 0\\
\end{array}\right), 
B=(e_{2},e_{3},-)=\left(\begin{array}{ccc}
\alpha & \gamma & \mu\\
\beta & \delta & \nu\\
0 & 0 & 0\\
\end{array}\right),
\end{displaymath}
\begin{displaymath}
C=(e_{3},e_{1},-)=\left(\begin{array}{ccc}
\kappa &-x -\alpha & \xi\\
\chi & -y- \beta & \beta\\
0 & 0 & 0\\
\end{array}\right).
\end{displaymath}

The verification of the defined relation of  LTS, leads to the
contradiction of the condition that $ dim \mathfrak{M'}=2 $.

Let $ \mathfrak{M'}=<e_{1},e_{2}> $-be a two-dimensional Abelian ideal and
$ e_3 $- the vector completing $ e_1 $, $e_2 $ up to the basis. Then:

\begin{displaymath}
A=(e_{1},e_{2},-)=\left(\begin{array}{ccc}
0 & 0 & a\\
0 & 0 & b\\
0 & 0 & 0\\
\end{array}\right), 
B=(e_{2},e_{3},-)=\left(\begin{array}{ccc}
\alpha & \gamma & \mu\\
\beta & \delta & \nu\\
0 & 0 & 0\\
\end{array}\right),
\end{displaymath}
\begin{displaymath}
C=(e_{3},e_{1},-)=\left(\begin{array}{ccc}
\kappa &-a -\alpha & \xi\\
\chi & -b- \beta & \beta\\
0 & 0 & 0\\
\end{array}\right).
\end{displaymath}
 
 Deforming the vector $ e_1 $ in the limit of the subspace $ <e_{1},e_{2}>$,
the matrix $ A $ can be reduced to the form $ a=b=0 $ or $ a=1 $, $ b=0 $.

The verification of the defined relation of the  LTS, in the second
 case leads to the following realization of the operators $ A $, $ B $, $ C $:

\begin{displaymath}
A=0, 
B=(e_{2},e_{3},-)=\left(\begin{array}{ccc}
0 & 0 & 1\\
0 & 0 & 0\\
0 & 0 & 0\\
\end{array}\right),
\end{displaymath}
\begin{displaymath}
C=(e_{3},e_{1},-)=\left(\begin{array}{ccc}
0 & 0 & 0\\
0 & 0 & 1\\
0 & 0 & 0\\
\end{array}\right).
\end{displaymath} \; \; \; \; \; \; \; \; \; \; \; \; \; ( type V)

\begin{displaymath}
A=0, 
B=(e_{2},e_{3},-)=\left(\begin{array}{ccc}
0 & 1 & 0\\
0 & 0 & \pm 1\\
0 & 0 & 0\\
\end{array}\right),\;
C=0.
\end{displaymath}\; \; \; \; \; \; \; \; \; \; \; \;  (type VI)
\end{enumerate}
In conclusion to the conducted examination we have the following theorem:

\underline{{\bf Theorem 4.1.}} Let $ \mathfrak{M}=<e_{1},e_{2},e_{3}> $- be 
a solvable
 LTS of dimension 3, $ \mathfrak{G}$- its canonical enveloping Lie
algebra(solvable), and let $ A, B, C:\mathfrak{M} \longrightarrow \mathfrak{M}$
the linear transformations of the view: $ A=(e_{1},e_{2},-) $,  $ A=(e_{1},e_{2},-) $,  $ B=(e_{2},e_{3},-) $,  $ C=(e_{3},e_{1},-) $: with isomorphism 
accuracy, one can find the possibility of the following types:
\begin{itemize}
\item Type I. $ \mathfrak{M}$- Abelian Lie triple system.

\item Type II. \begin{displaymath}
A=0, 
C=0,
B=(e_{2},e_{3},-)=\left(\begin{array}{ccc}
0 & 0 & 1\\
0 & 0 & 0\\
0 & 0 & 0\\
\end{array}\right)
\end{displaymath} 
$ \mathfrak{G}=<e_{1}, e_{2}, e_{3}, e_{4}> $- four-dimensional 
non-decomposable nilpotent Lie algebra with defined relations
$$
[e_{2},e_{3}]=e_{4}, [e_{3},e_{4}]=-e_{1}
$$
(this is $ g_{4,1} $ algebra in Mubaraczyanov classification\cite{mub2}).
\item Type III. $ \mathfrak{M} $ is a direct product of a two-dimensional 
solvable  LTS $ <e_{1},e_{2}> $, and a one-dimensional Abelian 
$ <e_{3}> $ :
 
\begin{displaymath}
A=(e_{1},e_{2},-)=\left(\begin{array}{ccc}
0 & \pm 1 & 1\\
0 & 0 & b\\
0 & 0 & 0\\
\end{array}\right), 
B=0,
C=0
\end{displaymath}

$ \mathfrak{G}=<e_{1}, e_{2}, e_{3}, e_{4}> $  four-dimensional solvable and
decomposable Lie algebra, with defined relations:
$$
[e_{1},e_{2}]=e_{4}, [e_{2}, e_{4}]= \pm e_{1}
$$
moreover $ \mathfrak{G}=<e_{1},e_{2}, e_{4}> \oplus <e_{3}> $, where 
$ <e_{1}, e_{2}, e_{4}> $- three-dimensional solvable Lie (algebra
 $ g_{3,4 \setminus 5}$ in Mubaraczyanov classification \cite{mub2}).
\item Type IV. \begin{displaymath}
A=(e_{1},e_{2},-)=\left(\begin{array}{ccc}
0 & \pm 1 & 1\\
0 & 0 & 0\\
0 & 0 & 0\\
\end{array}\right), 
B=0,
C=(e_{3},e_{1},-)=\left(\begin{array}{ccc}
0 & -1 & \pm 1\\
0 & 0 & 0\\
0 & 0 & 0\\
\end{array}\right)
\end{displaymath},
 
$ \mathfrak{G}=<e_{1}, e_{2}, e_{3}, e_{4}> $- four-dimensional solvable and
non-decomposable Lie algebra, with defined relations:
$$
[e_{1},e_{2}]=e_{4}, [e_{2}, e_{4}]= \pm e_{1}
$$
$$
[e_{1},e_{3}]=\pm e_{4}, [e_{3}, e_{4}]= - e_{1}
$$
(algebra $ g_{4,5 \setminus 6}$ in Mubaraczyanov classification \cite{mub2}).

\item Type V  
\begin{displaymath}
B=\left(\begin{array}{ccc}
0 &  1 & 0\\
0 & 0 & \pm 1\\
0 & 0 & 0\\
\end{array}\right), 
A=C=0
\end{displaymath}

$ \mathfrak{G}=<e_{1}, e_{2}, e_{3}, e_{4}> $- four-dimensional solvable 
non-decomposable Lie algebra with defined relations:
$$
[e_{2},e_{3}]=e_{4}, [e_{2}, e_{4}]= - e_{1}
$$
$$
[e_{3},e_{4}]=\mp e_{2}
$$
(algebra $ g_{8 \setminus 9}$ in Mubaraczyanov classification \cite{mub2}).

\item Type VI 
\begin{displaymath}
A=0, 
B=(e_{2},e_{3},-)=\left(\begin{array}{ccc}
0 & 0 & 1\\
0 & 0 & 0\\
0 & 0 & 0\\
\end{array}\right),
C=(e_{3},e_{1}, -)=\left(\begin{array}{ccc}
0 & 0 & 0\\
0 & 0 & 1\\
0 & 0 & 0\\
\end{array}\right)
\end{displaymath}

$ \mathfrak{G}=<e_{1}, e_{2}, e_{3}, e_{4},e_{5}> $- five-dimensional solvable 
non-decomposable Lie algebra, with defined relations:
$$
[e_{1}, e_{2}]=e_{4}, [e_{1}, e_{3}]= - e_{5}
$$
$$
[e_{3}, e_{4}]=- e_{1}, [e_{3}, e_{5}]=-e_{2}
$$
(as a result we obtain an extension of four-dimensional Abelian ideal 
 $ \mathfrak{G}=< e_{1}, e_{2}, e_{4}, e_{5}>$  by means of $ <e_{3}>$,
algebra $ g_{4, 13} $ in Mubaraczyanov classification \cite{mub2}).
\item Type VII.\begin{displaymath}
A=0, 
B=(e_{2},e_{3},-)=\left(\begin{array}{ccc}
1 & 0 & 0\\
0 & 0 & 0\\
0 & 0 & 0\\
\end{array}\right),
C=(e_{3},e_{1}, -)=\left(\begin{array}{ccc}
0 & -1 & 0\\
0 & 0 & 0\\
0 & 0 & 0\\
\end{array}\right)
\end{displaymath} 

\end{itemize}
  
Lie algebra $ \mathfrak{G}=<e_{1}, e_{2}, e_{3}, e_{4},e_{5}> $- five-
dimensional solvable non-decomposable Lie algebra, with defined relations:
$$
[e_{2}, e_{3}]=e_{4}, [e_{1}, e_{3}]=  e_{5}
$$
$$
[e_{1}, e_{4}]=- e_{1}, [e_{2}, e_{5}]=-e_{1}, [e_{4},e_{5}]=e_{5}
$$
(algebra $ g_{4, 11} $ in Mubaraczyanov classification \cite{mub3,mub4}).

\section{CLASSIFICATION OF SPLITTING 3-DIMENSIONAL LIE TRIPLE SYSTEMS }

 Let $ \mathfrak{M}=\mathfrak{M_{1}}\dotplus \mathfrak{M_{2}} $- be 
a splitting 3-dimensional 
 LTS, where  $ \mathfrak{M_{1}} \cong \mathbb{R}$- is a one dimensional solvable ideal in  $ \mathfrak{M}$ and  $ \mathfrak{M_{2}}$ be a 2-dimensional simple 
 LTS. Introduce in consideration a basis 
$(e_{1}, e_{2}, e_{3})$ in  $ \mathfrak{M}$ such that   $ \mathfrak{M_{1}}=<e_1> $ 
and  $ \mathfrak{M_{2}}=<e_{2}, e_{3}>$ and linear operators
    $ A, B, C:\mathfrak{M} \longrightarrow \mathfrak{M}$
 such that $ A=(e_{1},e_{2},-) $,  $ A=(e_{1},e_{2},-) $,  $ B=(e_{2},e_{3},-) $,  $ C=(e_{3},e_{1},-) $ using the process apply in the previous case one can
obtain the following theorem:
\\

\underline{{\bf Theorem 5.1.}} The following situation are possible and
non isomorphic:

\begin{itemize}
\item {\bf Type 1.} $ \mathfrak{M}=\mathbb{R}\oplus \mathfrak{M_{2}}$- direct sum of one dimensional Abelian ideal and
2-dimensional simple ideal in  $ \mathfrak{M}$
where $ \mathfrak{M_{2}}$ is a simple a simple 2-dimensional  LTS 
of the view $so(3)/so(2), sl(2, \mathbb{R})/so(2), sl(2, \mathbb{R})/\mathbb{R}$

\item {\bf Type 2.} \begin{displaymath}
A=\left(\begin{array}{ccc}
0 & -1 & 0\\
0 & 0 & 0\\
0 & 0 & 0\\
\end{array}\right), 
B=\left(\begin{array}{ccc}
0 & 0 & 0\\
0 & 0 & -1\\
0 & 1 & 0\\
\end{array}\right),
C=\left(\begin{array}{ccc}
0 & 0 & 1\\
0 & 0 & 0\\
0 & 0 & 0\\
\end{array}\right)
\end{displaymath} 
\\

$\mathfrak{M_{2}}$ is a simple 2-dimensional  LTS of the view
$so(3)/so(2)$

\item {\bf Type 3.}
 
\begin{displaymath}
A=\left(\begin{array}{ccc}
0 &  1 & 0\\
0 & 0 & 0\\
0 & 0 & 0\\
\end{array}\right), 
B=\left(\begin{array}{ccc}
0 & 0 & 0\\
0 & 0 & 1\\
0 & -1 & 0\\
\end{array}\right),
C=\left(\begin{array}{ccc}
0 & 0 & -1\\
0 & 0 & 0\\
0 & 0 & 0\\
\end{array}\right)
\end{displaymath}

$\mathfrak{M_{2}}$ is a simple 2-dimensional  LTS of the view
$sl(2, \mathbb{R})/so(2)$

\item {\bf Type 4}. \begin{displaymath}
A=\left(\begin{array}{ccc}
0 & - 1 & 0\\
0 & 0 & 0\\
0 & 0 & 0\\
\end{array}\right), 
B=\left(\begin{array}{ccc}
0 & 0 & 0\\
0 & 0 & 1\\
0 & 1 & 0\\
\end{array}\right),
C=\left(\begin{array}{ccc}
0 & 0 & - 1\\
0 & 0 & 0\\
0 & 0 & 0\\
\end{array}\right)
\end{displaymath},
 
$\mathfrak{M_{2}}$ is a simple 2-dimensional  LTS of the view
$sl(2, \mathbb{R})/so(2)$

\item {\bf Type 5.}  
\begin{displaymath}
A=-C=\left(\begin{array}{ccc}
0 & - \frac{1}{4} & \frac{1}{4}\\
0 & 0 & 0\\
0 & 0 & 0\\
\end{array}\right), 
B=\left(\begin{array}{ccc}
-\frac{1}{2} & 0 & 0\\
0 & 0 & 1\\
0 & 1 & 0\\
\end{array}\right)
\end{displaymath}
\\
$\mathfrak{M_{2}}$ is a simple 2-dimensional  LTS of the view
$sl(2, \mathbb{R})/so(2)$

\item {\bf Type 6.} 
\begin{displaymath} 
A=C=\left(\begin{array}{ccc}
0 & - \frac{1}{4} & -\frac{1}{4}\\
0 & 0 & 0\\
0 & 0 & 0\\
\end{array}\right), 
B=\left(\begin{array}{ccc}
\frac{1}{2} & 0 & 0\\
0 & 0 & 1\\
0 & 1 & 0\\
\end{array}\right)
\end{displaymath} 
$\mathfrak{M_{2}}$ is a simple 2-dimensional  LTS of the view
$sl(2, \mathbb{R})/so(2)$
\end{itemize}
  
The proof is somewhat intricate calculation as done in the section above.

{\bf Aknoledgment:} This paper was able to be achieved, thanks to the scholarship obtain from the Agence Universitaire de la Francophonie.

\end{document}